\documentclass{article}
\usepackage{amsfonts,amssymb,latexsym,amsmath}
\newcounter{num}[section]
\newcommand{\Num}{\refstepcounter{num}%
\textbf{\arabic{section}.\arabic{num}}}

\newcommand{\Theorem}{\textbf{Theorem~}}
\newcommand{\Proof}{\textbf{Proof}}
\newcommand{\Def}{\textbf{Definition~}}

\newcommand{\Lemma}{ \textbf{Lemma~}}
\newcommand{\Ex}{  \textbf{Example~}}

\newcommand{\Prop}{\textbf{Proposition~}}
\newcommand{\Cor}{ \textbf{ Corollary~}}

\newcommand{\Ax}{{\mathfrak A}}
\newcommand{\Bx}{{\mathfrak B}}

\newcommand{\Kc}{{\cal K}}
\newcommand{\Oc}{{\cal O}}
\newcommand{\Xc}{{\cal X}}
\newcommand{\Fc}{{\cal F}}
\newcommand{\Ac}{{\cal A}}
\newcommand{\Bc}{{\cal B}}

\newcommand{\Ch}{{{\mathfrak C}{\mathfrak h}}}
\newcommand{\al}{{\alpha}}

\newcommand{\la}{{\lambda}}

\newcommand{\tG}{{\widetilde{G}}}

\newcommand{\eps}{{\varepsilon}}

\newcommand{\UTn}{{\mathrm{UT}}(n,\Fq)}
\newcommand{\Tn}{{\mathrm{T}}(n,\Fq)}
\newcommand{\Irr}{{\mathrm{Irr}}}

\newcommand{\Fq}{{\Bbb F}_q}
\newcommand{\Cb}{{\Bbb C}}
\newcommand{\Zb}{{\Bbb Z}}

\renewcommand{\leq}{\leqslant}
\renewcommand{\geq}{\geqslant}

\newcommand{\Res}{{\mathrm{Res}}}
\newcommand{\Ind}{{\mathrm{Ind}}}
\newcommand{\SInd}{{\mathrm{SInd}}}
\newcommand{\rt}{{\mathrm{right}}}
\newcommand{\lt}{{\mathrm{left}}}

\begin{document}
\Large

\title{Restriction and induction for supercharacters  of finite groups of triangular type}
\author{A.N.Panov
\thanks{The paper is produced by  financial support of the Ministry of education and science of the Russian Federation within the basic part of the state task (project No. 204). The work is supported by RFBR grants РФФИ 16-01-00154-a and 14-01-97017-Volga-region-a.}}
\date{}

 \maketitle

\begin{abstract}
It is proved the  restriction of any  supercharacter of a finite group of triangular type on its subgroup is a sum of supercharacters with nonnegative integer coefficients. We define a superinduction and prove the analog of  Frobenius reciprocity formula for supercharacters.
\end{abstract}

\section{Introduction}
The problem of classification of irreducible characters (representations) is the main problem in the representation theory of finite
groups. It is well known that a character is constant on  conjugacy  classes, the number of irreducible characters is equal to the number of conjugacy classes. The values of irreducible characters on conjugacy classes form the quadratic table which is called a character table.

However, for some finite groups, like  the unitriangular group
 $\mathrm{UT}(n,\Fq)$, the problem of classification of irreducible characters is a very complicated, "wild"\, problem. In the paper,  \cite{DI} ~P.Diaconis and I.Isaaks proposed the concept of a supercharacter theory.  Roughly speeking a supercharacter theory is a
 first approximation of the theory of irreducible characters.  Together with supercharacters there are defined superclasses;
       the values of supercharacters are constant on superclasses and form the quadratic table, which is called a supercharacter table.  Originally, any finite group have several supercharacter theories;  preference is given to the supercharacter theory that produces better approximation of the theory of irreducible characters.

Historically, the first nontrivial example of a supercharacter theory was the theory of basic characters of   $\UTn$ constructed in the series of papers of K.Andre  \cite{A1,A2,A3,A4} (see also \cite{Yan}).
    In the paper  \cite{DI},  ~P.Diaconis and I.Isaaks constructed the supercharacter theory for algebra groups   (see example  \ref{extree}) which generalizes the Andre theory.
 Many papers were devoted to investigating different supercharacter theories. Highlight  a few of them:  supercharacter theory for abelian groups with application to the number theory \cite{Number-1, Number-2},
superinduction for algebra groups  \cite{ThiemB, ThiemRest, MT}, application to the problem of random walking on groups  \cite{Walk}, supercharacter theory for semidirect products  \cite{H}. There is  the bibliography  on this topic in the paper \cite{VERY}.

 Notice that for algebra groups the restriction of a supercharacter  on its algebra subgroup  is a sum of supercharacters
 with nonnegative integer coefficients. At the same time, in general, the induced character from a supercharacter is not a sum of supercharacters of group (even with complex coefficients). But if you change induction by superinduction, then the constructed   character will be the sum of supercharacters (with nonnegative rational coefficients), and   there is an analog of  Frobenius reciprocity theorem (see \cite{DI}).

  In the papers  \cite{P1,P2}, there were constructed a supercharacter theory for finite groups of triangular type that generalizes
  the theory of  P.Diaconis and I.Isaaks.  A goal of this paper is to transfer the mentioned statements on restriction and superinduction for algebra groups to the case of finite groups of triangular type. The main results are formulated in theorems
\ref{theoremrestr} and \ref{last}.

\section{Overview of supercharacter theory}

Let $G$ be a group, ~$1\in G$ be the unit element, ~ $\Irr(G)$ be the set of irreducible characters (representations) of the group ы $G$.
   Suppose that we have two partitions
    \begin{equation}\label{partitionIrrG}
 \Irr(G) = X_1\cup \cdots \cup X_m,\quad X_i\cap X_j=\varnothing,
 \end{equation}
  \begin{equation}\label{partitionG}
G = K_1\cup \cdots \cup K_m,\quad K_i\cap K_j=\varnothing.
\end{equation}
 Notice that the number of components of these partitions are common. To each $X_i$ we correspond  the character of  group $G$ by the formula
 \begin{equation}\label{sigma}
   \sigma_{i} = \sum_{\psi\in X_i} \psi(1) \psi.
 \end{equation}
 \Def\Num\label{twopart}. Two partitions  $\Xc=\{X_i\}$ и $\Kc=\{K_j\}$ is said to define  a supercharacter theory for  $G$
 if each character  $\sigma_i$ is constant on each  $K_j$. In this case,   $\{ \sigma_i\}$  are called  \emph{supercharacters}, and $\{K_j\}$ -- \emph{superclasses}.  The table of values $\{ \sigma_i(K_j)\}$ is called a supercharacter table.

 The following statement is important for constructing supercharacter theories.\\
\Prop\Num\label{dilemma}~ \cite[Lemma 2.1]{DI}. Suppose that we have a system of disjoint characters   $\Ch=\{\chi_1,\ldots,\chi_m\}$ and a partition   $\Kc=\{K_1,\ldots, K_m\}$ of the group $G$. Assume that each character  $\chi_i$ is constant on each  $K_j$.  Denote by  $X_i$ the \emph{ support of the character } $\chi_i$ (i.e., the set of all irreducible components of $\chi_i$).   Then the following conditions are equivalent:\\
1) $\{1\}\in \Kc$,\\
2)  the system of subsets  $\Xc=\{X_i\}$  form a partition of  $\Irr(G)$; two partitions  $\Xc$ and  $\Kc$ defines a supercharacter theory of the group  $G$. Moreover,  each $\chi_i$ is up to a constant factor equal to  $\sigma_i$.
Simplifying language, we refer to  $\chi_i$ as \emph{supercharacters}.\\
\\
\Ex\Num. The system of irreducible characters  $\{\chi_i\}$ and the system of conjugacy classes
$\{K_i\}$.\\
\Ex\Num. Two supercharacters  $\chi_1=1_G$, ~ $\chi_2=\rho-1_G$ (here $\rho$ is a character of the regular representation) and two superclasses
$K_1=\{1\}$, ~ $K_2=G\setminus \{1\}$. The supercharacter table has the form
\begin{center}\begin{tabular}{|c|c|c|}
                         \hline
                          & $\chi_1$ &  $\chi_2$\\
                          \hline
      $K_1$ & 1  &  $|G|-1$\\
      \hline
                $K_2$ & 1  &  -1\\
              \hline
               \end{tabular}
\end{center}
\Ex \Num\label{extree}. The supercharacter theory for algebra groups  \cite{DI} (see also  \cite{P3}).
By definition,  \emph{an algebra group} is a group  $G=1+J$, where  $J$ is an associative finite dimensional nilpotent algebra over the finite field  $k=\Fq$.    The superclass of the element  $1+x$ is defined as $1+\omega$, where $\omega $ is a left-right  $G\times G$-orbit of the element  $x\in J$.
There are also right and left actions of the group  $G$ on the dual space  $J^*$ defined by the formulas  $\la g(x) = \la(gx)$ and $g\la(x)=\la(xg)$.
Let  $G_{\la,\rt}$ be the stabilizer of  $\la\in J^*$ with respect to the right action of $G$ on $J^*$.

Fix a nontrivial character  $t\to \eps^t$ of the additive group  $\Fq$ with values in the multiplicative group  $\Cb^*$. The function $\xi :G_{\la,\rt}\to \Cb^*$ defined as $$\xi_\la(g) = \eps^{\la(g-1)}$$
is a linear character (one dimensional representation) of the group  $G_{\la,\rt}$.
 A \emph{supercharacter } of the algebra group  $G$ is the induced character
 \begin{equation}\label{chiag}
 \chi_\la = \Ind(\xi_\la, G_{\la,\rt}, G).
\end{equation}
  The systems of supercharacters  $\{\chi_\la\}$ and superclasses  $\{1+GxG\}$, where $\la$  and $x$  run through the systems of double  $G\times G$-orbits in  $J^*$ and $J$ respectively, define a supercharacter theory for the group  $G$.

Let us observe the problems of restriction and induction in the theory supercharacters for algebra groups.
Let  $G =1+J$ be an algebra group. If  $J'$ is an arbitrary subalgebra in  $J$, then the subgroup  $G'=1+J'$ is called an \emph{ algebra subgroup} in $G$. \\
It is proved in the paper  \cite{DI} that the restriction of  supercharacter  $\chi_\la$ on the algebra subgroup  $G'$ is a sum of supercharacters of the subgroup  $G'$ with nonnegative integer coefficients.

  Suppose that $\phi$ is a superclass function on  $G'$ (i.e., the function constant on superclasses of  $G'$). Extend  $\phi$ to the function  $\dot{\phi}$ on $G$ letting it equal to zero outer $G'$.  \\
By definition, a \emph{superinduction} of $\phi$ is a function $\SInd\,\phi$   on the group $G$ defined as follows
$$ \SInd\,\phi(1+x) = \frac{1}{|G|\cdot|G'|}\sum_{a,b\in G}\dot{\phi}(1+axb).$$
Easy to see that  $ \SInd\,\phi(1+x)$ is a superclass function on the group  $G$.
There is the standard scalar product on the group $G$ (and $G'$) defined as
$$(f_1,f_2) = \frac{1}{|G|}\sum_{g\in G} f_1(g)\overline{f_2(g)}.$$
The following theorem is  an supercharacter analog of  the Frobenius reciprocity theorem for algebra groups.\\
\Theorem\Num\,\cite{DI}. Let $\phi$ be a superclass function on $G'$, and $\psi$  be a superclass function on  $G$. Then
$(\SInd\,\phi,\psi) = (\phi, \Res\,\psi).$

\section{Supercharacters of finite groups of triangular type}

In this section, we highlight the supercharacter theory for finite groups of triangular type constructed by the author in the papers \cite{P1,P2}.

 Let $H$ be a group, and  $J$ is an associative algebra over a field  $k$. Suppose that there defined the left  $h,x\to hx$ and the right  $h,x\to xh$ linear actions of the group  $H$ on $J$. Assume that for every  $h\in H$ and $x,y\in J$ the following conditions are fulfilled:
\begin{enumerate}
\item $h(xy)=(hx)y$ и $(xy)h=x(yh)$,
\item $x(hy)=(xh)y$.
\end{enumerate}
On the set  $$G=H+J=\{h+x:~ h\in H, ~ x\in J\}$$
we define an operation of multiplication
\begin{equation}\label{mult}
g_1g_2=(h_1+x_1)(h_2+x_2)=h_1h_2+ h_1x_2+x_1h_2+x_1x_2.
\end{equation}
If  $J$ is a nilpotent algebra over the  field  $k$, then $G$ is a group with respect to operation  (\ref{mult}). If the group  $H$ is finite, the field $k$ is finite, and  $J$ is finite dimensional algebra over the field  $k$,  then the group  $G$ is finite. \\
\Def\Num. Under all these conditions we call the group  $G$  a \emph{finite group of triangular type} if the group  $H$ is abelian
and  $\mathrm{char}\, k$ does not divide  $|H|$.

Let  $G=H+J$ be a finite group of triangular type. The group algebra $kH$  is commutative and  semisimple
(according to Maschke's theorem). Therefore, $kH$  is a sum of fields.
There exists a system of primitive idempotents  $\{e_1,\ldots, e_n\}$ such that   \begin{equation}\label{gralg} kH=k_1e_1\oplus\ldots \oplus k_ne_n, \end{equation} where $k_1,\ldots,k_n$ are  field extensions of  $k$. Any idempotent in  $kH$ is a sum of primitive idempotents.

The direct sum  $A=kH\oplus J$ has an algebra structure  with respect to the multiplication  (\ref{mult}). The group  $G$ is a subgroup of the group  $A^*$ of invertible elements of $A$, which is considered in the Example \ref{exthree}.
 Observe that the group  $G$  decomposes into a product $G=HN$ of the subgroup  $H$ and the normal subgroup   $N=1+J$ that is an algebra group.\\
\Ex\Num. Algebra group $G=1+J$.\\
\Ex\Num. $G=\left\{\left(\begin{array}{cc}a&b\\
0&1\end{array}\right):~~ a,b\in \Fq,~ a\ne 0\right\}$.\\
\Ex\Num\label{exthree}. Let $A$ be an associative unital finite dimensional algebra  of reduced type over the finite field   $\Fq$ of $q$ elements \cite[\S 6.6]{Pi}. By definition, the  algebra $A$ is reduced if  its factor algebra with respect to the radical    $J=J(A)$ is a direct sum of division algebras. According to Wedderburn's theorem  \cite[\S 13.6]{Pi}, any division algebra over a finite field is commutative. Then the algebra  $A/J$ is commutative.
There exists a semisimple subalgebra $S$ such that
$A = S\oplus J$ (see \cite[\S 11.6]{Pi}). In our case,  $S$ is commutative.
The group  $G=A^*$  on the invertible elements of  $A$ is a finite group of triangular type    $G=H+J$, where  $H=S^*$. If  $A$ is the algebra of triangular matrices, then $G=\Tn$ is the triangular group.

 Let $G=H+J$ be a finite group of triangular type. Consider the group $\tG$ of triples  $\tau=(t,a,b)$, where $t\in H$,~ $a,b\in N$, with operation  $$(t_1,a_1,b_1)\cdot (t_2,a_2,b_2) = (t_1t_2,~ t_2^{-1}a_1t_2a_2, ~ t_2^{-1}b_1t_2b_2).$$
The group   $\tG$ acts on   $J$  as follows
$$ \rho_\tau(x) = taxb^{-1}t^{-1}.$$
In the dual space  $J^*$ a representation of the group  $\tG$ is defined as usual  $$\rho^*_\tau\la(x) = \la(\rho(\tau^{-1})(x)).$$   In the space $J^*$ there are also left and right linear actions of the group  $G $ by the formulas  $b\la(x)= \la(xb)$ and $\la a(x) = \la(ax)$. Then  $\rho_\tau(\la) =
tb\la  a^{-1}t^{-1}$.

 For any idempotent  $e\in kH$ we denote by  $A_e$ the subalgebra  $eAe$.
The subalgebra  $J_e=eJe\subset J$ is a radical in  $A_e$. Denote  $e'=1-e$. There is a the Pierce decomposition
 $$J= eJe\oplus eJe'\oplus e'Je\oplus e'Je'.$$
 The dual space $J_e^*$ is identified with the subspace in $J^*$ that consists of all linear forms equal to zero on all components of the  Pierce decomposition except the first one.
 Observe that since the group   $H$ is abelian, we have  $he=eh=ehe$ for all  $h\in H$. The subset $H_e=eHe$ is a subgroup in the algebra of invertible elements of the algebra  $A_e$. The group  $G_e=eGe=H_e+J_e$ is a group of triangular type; and it is associated with the algebra  $A_e$ in the same way as    $G$ is associated with  $A$.
 The map $h\to he$ is a homomorphism of the group  $H$ onto  $H_e$, its kernel is the subgroup
\begin{equation}\label{he} H(e) =\{h\in H: ~ he=e\}.\end{equation}

The following definition and the one of \cite{P1,P2} differs in form, but equivalent.
\\
 \Def\Num\label{sing}. We say that  a $\rho_{\tG}$-orbit $\Oc$ is  \emph{singular} (with respect to  $H$) if  $\Oc\cap J_e\ne\varnothing$ for some idempotent  $e\ne 1$ in $kH$. Otherwise, the orbit  $\Oc$ is called  \emph{regular} (with respect to  $H$).
   Elements of singular (regular) orbits are called singular (regular).  Similarly defined singular and regular orbits and elements in $J^*$.

  The subgroup  $H(e)$ admits the following characterization.\\
\Prop\Num\label{charac}~\cite[Lemma 2.5]{P2}.
1)~  $H(e)   = H_{y,\mathrm{right}}\cap H_{y,\mathrm{left}}$ for any regular (with respect to  $H_e$) element  $y\in J_e$.\\
2)~ $H(e) = H_{\la,\mathrm{right}}\cap H_{\la,\mathrm{left}}$ for any regular (with respect to $H_e$) element  $\la\in J_e^*$.

It is proved in the papers  \cite{P1,P2} that for any  $\tG$-orbit $\Oc$ in $J$ the following statements are true. \\
 \Prop\Num\label{state}. 1) The intersection  $\Oc\cap J_e$ is an  $\tG_e$-orbit in $J_e$ for any idempotent  $e\in kH$.\\
 2) There exists a unique idempotent  $e\in kH$ such that  $\Oc\cap J_e$ is a regular  $\tG_e$-orbit in  $J_e$ (with respect to  $H_e$). \\
 Similar statements are true for  $\tG$-orbits in $J^*$.

Well known that  for any finite transformation group of finite dimensional linear space   $V$ defined over a finite field, the numbers of orbits in  $V$ and  $V^*$ are equal \cite[Lemma 4.1]{DI}. By this statement and above properties of $\tG$-orbits,
the numbers of regular (singular)  $\tG$-orbits in  $J$ and $J^*$ are equal \cite[Proposition 2.11]{P1}.

Turn to a definition of superclasses in the group  $G$. For each  $g\in G$ and $(t,a,n)\in\tG$ consider the element \begin{equation}\label{rtau}
R_\tau(g) = 1+ ta(g-1)b^{-1}t^{-1}
\end{equation}
in the algebra  $A=kH+J$.
  If $g=h+x$, then
$ R_\tau(g)=h\bmod J$. Hence  $ R_\tau(g)\in G$.
The formula  (\ref{rtau}) determines the action of the group  $\tG$ on $G$. \\
\Def\Num. We say that    $\tG$-orbits in $G$ are  \emph{superclasses}.

 The group  $G$  splits into superclasses.
Denote by $\Bx$ the set of triples  $\beta = (e,h, \omega)$, where
 $e$ is an idempotemt in   $kH$, ~ $h\in H(e)$,  and  $\omega$ is a regular (with respect to  $H_e$) ~ $\tG_e$-orbit in  $J_e$. All element  $h+\omega$ belong to a common superclass  \cite[следствие 3.2]{P2} denoted by  $K_\beta$.\\
\Theorem\Num~ \cite[теорема 3.3]{P2}. The correspondence  $\beta\to K_\beta$ is a bijection of the set of triples  $\Bx$ onto the set of superclasses in $G$.

Denote by $\Ax$ the set of triples  $\al = (e,\theta, \omega^*)$, where
 $e$ is an idempotent in  $kH$, ~ $\theta $ is a linear character (one dimensional representation) of the subgroup  $ H(e)$,  and $\omega^*$ is a regular (with respect to $H(e)$) ~ $\tG_e$-orbit in $J^*_e$.
 Since the subgroup  $H(e)$ is abelian, the number of its linear characters is equal to the number of its elements. The number of regular (with respect to  $H_e$)  ~ $\tG_e$-orbits in  $J_e$ and  $J_e^*$ are equal. Therefore  $|\Ax|= |\Bx|$.

Turn to a construction of supercharacters.
Let  $\al = (e,\theta, \omega^*)\in \Ax$, choose $\la\in \omega^*$.
 Consider the subgroup   $G_\al = H(e)\cdot N_{\la,\mathrm{right}}$, where
 $N_{\la,\mathrm{right}}$ is a stabilizer for $\la$ of the right action of  $N=1+J$ on $J^*$.
 Any element  $g\in G_\al$ can be presented in the form   $g=h+x$,  where  $h\la=\la h=\la$ and $\la (xJ)=0$. Hence
 \begin{equation}\label{Gal}
   G_\al = H(e) + J_{\la,\rt}.
 \end{equation}

Fix a nontrivial character  $t\to  \varepsilon^{t}$  of the additive group  $\Fq$with values in multiplicative group  $\Cb^*$.
By given  a triple  $\al = (e,\theta, \omega^*)$ and $\la \in\omega^*$, define   a linear character of the subgroup  $G_\al$ by the formula
\begin{equation}\label{thela}
    \xi_{\theta,\la}(g) = \theta(h)\varepsilon^{\la(x)},
\end{equation}
where $g=h+x$,~ $h\in H(e)$ и $x\in J_{\la,\mathrm{right}}$.

The induced character
\begin{equation}\label{Indchi}
\chi_\al = \Ind(\xi_{\theta,\la}, G_\al, G)
\end{equation}
is called a  {\it supercharacter}.\\
\Theorem\Num ~\cite[Propositions 4.1, 4.4]{P1,P2}. \\
1) Supercharacters $\{\chi_\al:~ \al \in\Ax\}$ are pairwise disjoint;\\
 2) Each supercharacter   $\chi_\al$ is constant on each superclass $K_\beta$;\\
 3) $\{1\}$ is a superclass  $K(g)$ for $g=1$.

 Applying Proposition  w \ref{dilemma}, we conclude.\\
 \Theorem\Num ~\cite[Theorem 4.5]{P2}. The systems of supercharacters  $ \{\chi_\al|~ \al\in \Ax\}$ and superclasses   $\{\Kc_\beta|~ \beta\in \Bx\}$  determines a supercharacter theory of the group    $G$.

\section{Restriction and superinduction for finite groups of triangular type}

Let  $G=H+J$ be a finite group of triangular type.
Let  $H$ be a subgroup of  $H$, ~ and  $J'$ be a subalgebra of $J$ invariant with respect to  left-right action of $H'\times H'$ on $J$. Then  $G'=H'+J'$ is a subgroup in  $G$, which will be called a \emph{subgroup of triangular tipe} in $G$. Our goal is to verify  that the restriction of a supercharacter on a subgroup of triangular type is  a sum of supercharacters of the subgroup with nonnegative integer coefficients.
\\
\Lemma\Num\label{andre}\,\cite[Lemma 6.1]{DI}. If $\Ac$ is a nilpotent associative algebra, $\Bc$ is a  subalgebra in $\Ac$, and $\Bc+\Ac^2=\Ac$, then $\Bc=\Ac$.\\
\Proof. For any positive integer  $n$ we have  $$\Ac^n=(\Bc+\Ac^2)^n\subseteq \Bc^n+\Ac^{n+1}\subseteq \Bc+\Ac^{n+1}.$$
Hence $\Bc+\Ac^n\subseteq \Bc+\Ac^{n+1}$. Since  $\Ac$ is a nilpotent algebra,   $\Ac^{n+1}=0$ for some  $n$. Then  $$\Ac=\Bc+\Ac^2\subseteq \Bc+\Ac^3 \subseteq\ldots \subseteq \Bc+\Ac^{n+1} =\Bc.$$
$\Box$\\
\Cor\Num\label{bhmax}. If   $J'$ is maximal $H\times H$-invariant subalgebra in $J$, then $J^2\subseteq J'$. In particular,  $J'$ is a ideal in  $J$.\\
\Proof.  Since   $J'\subseteq J'+J^2\subseteq J$, and   $J'$ is maximal  $H\times H$- invariant subalgebra in  $J$, then either
$J'=J'+J^2$, or $J'+J^2=J$. The second equality implies  $J'=J$, this contradicts to the assumption.
Therefore  $J'=J'+J^2$. Hence $J^2\subseteq J'$. $\Box$

Notice that the algebra $J$, as the algebra group  $N=1+J$, acts on $J^*$ by left and right transformations   $\la x(y)=\la(xy)$ and $x\la(y)=\la(yx)$. The set  $J\la$ is a left orbit of  $J$ on  $J^*$. \\
\Lemma\Num\label{lemann}~\cite[Lemma 4.2]{DI}.
$J_{\la,\rt}^\perp = J\la$ for any  $\la\in J^*$.
\\
\Proof. If $x\in J_{\la,\rt}$, then $\la(xJ)=0$. Hence $J\la(x)=0$. This verifies the inclusion $J\la\subset J_{\la,\rt}^\perp$.

To prove equality it is sufficient to show that the dimensions of subspaces
$J\la$ and $J_{\la,\rt}^\perp$ are equal. The dimension of the first one equals to  $\dim\,J-\dim\,J_{\la,\lt}$, and the second one  to $\dim\,J-\dim\,J_{\la,\rt}$. The coincidence of dimensions of left and right stabilizers in  $J$ can be verified  from the fact that the bilinear form  $\la(xy)$ provides a pairing of
$J/J_{\la,\lt}$ and $J/J_{\la,\rt}$. $\Box$

Let $\la\in J^*$ and  $H_0$ be a subgroup in  $H$ such that  $h_0\la=\la h_0=\la$ for every $h_0\in H_0$. Let $\theta_0$ be a linear character of $H_0$. Then $G_0=H_0+J_{\la,\rt}$ is a subgroup in  $G$, and the following formula
$$\xi_{{\theta_0}, \la}(g) ={\theta_0}(h_0)\eps^{\la(x)},\quad g=h_0+x$$
defines a linear character of the subgroup  $G_0$. Denote
$\chi_{{\theta_0},\la}=\Ind(\xi_{{\theta_0}, \la}, G_0,G)$.\\
\Lemma\Num\label{indind}. The induced character   $\chi_{{\theta_0},\la}$ is a sum of supercharacters.\\
\Proof. Let $f$ be the least idempotent in  $kH_0$ such that  $\la\in J_f^*$.
There exists the idempotent  $e\leq f$ such that the intersection of the  $\tG_f$-orbit of  $\la$ in $J_f^*$ is a regular
 $\tG_e$-orbit in  в $J_e^*$ (see Proposition  \ref{charac}).
Choose  $\la'$ in $J_e^*$ lying in the same  $\tG_f$-orbit of $\la$.
Then $\chi_{{\theta_0},\la}=\chi_{{\theta_0},\la'}$ (see \cite[Proposition 4.2]{P1,P2}).
The subgroup  $H(e)$ is contained in  $H_0$. Decompose  $\Ind({\theta_0},H_0,H(e))$ into a sum of linear characters
$$  \Ind({\theta_0},H_0,H(e)) = \sum_{i=1}^m \theta_i. $$
Then  $\chi_{{\theta_0},\la}$  is a sum of supercharacters constructed by  $(\bar{e},\theta_i,\omega'^*)$, where
$\omega'^*$  is an orbit of  $\la'$ with respect to  $\tG_e$. $\Box$

Let us state and prove the main theorem of this paper.\\
\Theorem\Num\label{theoremrestr}. The restriction of a supercharacter of the group  $G$ to its triangular type subgroup  $G'$ is a a sum of supercharacters of the subgroup  $G'$ with nonnegative integer coefficients.
\\
\Proof. There exists a chain of  $H'\times H'$-invariant subalgebras  $J'=J_1\subseteq \ldots\subseteq J_k=J$ such that
 each subalgebra  $J_i$ is a maximal  $H'\times H'$-invariant subalgebra in  $J_{i+1}$.
We complete the chain of subgroups  $G'=G_1\subseteq\ldots\subseteq G_k=H'+J$  by  $G_{k+1}=G=H+J$.
It is sufficient to verify the statement of the theorem in the case of restriction of supercharacter of the group $G_i$ to $G_{i-1}$; that is in the following two cases: 1)  $G=H+J$ and $G'=H'+J$, and 2) $G=H+J$  and $G'=H+J'$, where $J'$ is a maximal $H\times H$-invariant subalgebra in  $J$.

Let  $\chi_\al$ be a supercharacter of the group  $G=H+J$ constructed by the triple $\al=(e,\theta, \omega)$.\\
{\bf Case 1}.~ $G=H+J$ and $G'=H'+J$.  According to the well known theorem on restriction of  induced representation on  subgroup \cite[Theorem 44.2]{CR}, we have
\begin{equation}\label{restr}
\Res\,\chi_\al = \sum_s \Ind(\xi_{\theta,\la}^{(s)}, G^{(s)}_\al\cap G', G'),
\end{equation}
 where $ G^{(s)}_\al = sG_\al s ^{-1}$,~  $\xi_\al^{(s)}$ is a character of subgroup  $ G^{(s)}_\al\cap G'$ defied by the formula  $$\xi_{\theta,\la}^{(s)} (k) = \xi_{\theta,\la}(s^{-1}ks), $$ and $s$ runs through the set of representatives of double cosets   $G'\backslash G/G_\al$. It is sufficient to prove that each character included in the sum  (\ref{restr}) is a sum of supercharacters.

Under conditions of Case 1, we consider  that  $s$ belongs to  $H$ and runs through the set of representatives of the cosets of the subgroup  generated by   $H(e)$ and $H'$.
The subgroup  $G_{\al}^{(s)}\cap G'$ coincides with  $$ \left(H(e)\cap H'\right)+J_{s\la s^{-1},\rt}.$$
Choose  $\la\in\omega^*$. The character  $ \xi_{\theta,\la}^{(s)}$ of the subgroup  $G_{\al}^{(s)}\cap G'$ is calculated by the formula  $$ \xi_{\theta,\la}^{(s)}(h+x) = \theta(h)\varepsilon^{s\la s^{-1}(x)}.$$

 Let $\bar{e}$ be the least idempotent in  $kH'$ that  $\geq e$. Show that the subgroup  $ H(e)\cap H'$  coincides with   $$H'(\bar{e}) =\{h'\in H':~~ h'\bar{e} =\bar{e}\}.$$
From Proposition   \ref{charac} we see that  $H(e)\cap H'= H'_{\la,\mathrm{right}}\cap H'_{\la,\mathrm{left}}$.
Let  $\bar{\omega}$ be an orbit of  $\la$ in  $J^*_{\bar{e}}$ with respect to the group $\tG_{\bar{e}}$. Taking into account
Proposition  \ref{charac}, to prove of equality  $H'(\bar{e})=H(e)\cap H'$ it is sufficient to show that   $ \bar{\omega}$ is a regular with respect to  $H'_{\bar{e}}$. Assume the contrary. Let   $\bar{\omega}\cap J^*_{f}\ne\varnothing$ for some  idempotent  $f\in kH'$ and $f < \bar{e}$.
Recall that  $\la\in J^*_e$, and hence  $\bar{\omega}\cap J^*_{e}\ne\varnothing$.
Then  $\bar{\omega}\cap J^*_{ef}\ne\varnothing$ (see \cite[Lemma 2.5]{P1}).
Since  $ef\leq e$ and  the orbit  $\omega$ is regular in  $J^*_{e}$, we have $ef=e$. This implies  $e\leq f < \bar{e}$. It contradicts to minimality condition for  $\bar{e}$.
The equality  $H'(\bar{e})=H(e)\cap H'$ is proved.

 The orbits  $ s\omega^* s^{-1}$  (as the orbit  $ \omega^* $)  are regular in  $J_{\bar{e}}$ (with respect to   $H'_e$).
 Let $\bar{\theta}$ be the restriction of the character  $\theta$ on $H'(\bar{e})$. Then $\Ind(\xi_\al^{(s)}, G^{(s)}_\al\cap G', G')$ coincides with the supercharacter  $\chi_{\bar{\al}}$ of the subgroup  $G'$, where $\bar{\al}= (\bar{e}, \bar{\theta}, s\omega^* s^{-1})$. \\
 {\bf Case 2}. $G=H+J$  and $G'=H+J'$, where $J'$ is a maximal  $H\times H$-invariant subalgebra in  $J$.\\
The representation of   $H\times H$ in  $J$ is completely reducible. Since  $J^2\in J'$ (see Corollary \ref{bhmax}), every subspace of the form  $J'+W$, where $W$ in an invariant with respect to  $H\times H$, is a subalgebra (moreover, an ideal).
It follows from maximality of  $J'$ that
 $J=J'+W$, where $W$ is an irreducible with respect to  $H\times H$ subspace.
For  $J$  we have the  Pierce decomposition
\begin{equation}\label{pierce}
J = J_{11}\oplus J_{12}\oplus J_{21}\oplus J_{22},
\end{equation}
where $J_{11}= J_e = eJe$,~ $J_{12}=eJe'$,~ $J_{21} = e'Je$,~ $J_{22}= e'Je'$, ~ $e'=1-e$.
Similarly for  $J'$.
The Pierce components  are $H\times H$ invariant, therefore,  $W$ belongs to a unique the  Pierce component. \\
\emph{ Item 2.1.} Let  $W\subseteq J_{11}$. Then $J_{12}=J'_{12}$,~  $J_{21}=J'_{21}$,~ $J_{22}=J'_{22}$, and
$J_{11} = W\oplus J'_{11}$.

Apply decomposition  (\ref{restr}) to  $\Res\,\chi_\al$. Let us calculate  $G^{(s)}_\al\cap G'$.
Recall that  $G_\al=H(e)+J_{\la,\rt}$, where $H(e)=\{h\in H: he=eh=e\}$ and
$J_{\la,\rt}=\{y\in J:~ \la(yJ)=0\}$.

We can consider the representative   $s$ of $G'\backslash G/G_\al$ runs  through  $1+W\subset N_e$.
Hence for any  $h\in H(e)$ we obtain  $sh=(1+ewe)h=h+eweh=h+ewe=h+hewe=h(1+ewe)=hs$. Therefore $sH(e)s^{-1}=H(e)$.
Then  $$G^{(s)}_\al=sG_\al s^{-1} = H(e) + J_{s\la s^{-1},\rt}.$$
Denote  $\mu=s\la s^{-1}$. As   $\la$,  the element  $\mu$ belongs to  $J_e^*$.
 Likewise  (\ref{thela}), the subgroup   $G^{(s)}_\al = H(e)+J_{\mu,\rt}$ has the character  $\xi_{\theta,\mu}$.
We obtain
\begin{equation}\label{capcap}
G^{(s)}_\al\cap G'= H(e)+\Bc,
\end{equation}
where $\Bc= J_{\mu,\rt}\cap J'$.
  The character $\xi_\al^{(s)}$ of the subgroup  $H(e)+\Bc$  is calculated  $$\xi_\al^{(s)}(h+x) = \theta(h)\varepsilon^{\mu(x)} =\xi_{\theta,\mu}(h+x).$$
The formula  (\ref{restr}) has the form
\begin{equation}\label{restr1}
\Res\,\chi_\al = \sum_\mu \Ind(\xi_{\theta,\mu}, H(e)+\Bc, G'),
\end{equation}
where $\mu\in J_e^* $ runs through    $$\{s\la  s^{-1}:~~ s=1+W\}.$$
It is sufficient to prove that each summand in  (\ref{restr1}) is a sum of supercharacters of
 $G'$.

Denote by  $\mu'$ the natural projection of  $\mu$ on $J'^{*}$. As  $J^2\subset J'$,  we have
 $\Bc=\{y\in J':~ \mu'(yJ)=0\}$. The subalgebra  $\Bc$ is contained in the subalgebra
 $\Ac=J'_{\mu',\rt}= \{y\in J':~ \mu'(yJ')=0\}$.
 Since  $\mu\in J_e^*$ and $\mu'\in J_e'^*$, the subalgebras  $\Ac$ and $\Bc$ are invariant with respect to left and right multiplication by  $H(e)$ and graded with respect to the  Pierce decomposition (\ref{pierce}).

The subalgebras  $\Ac$ and  $\Bc$ may differ only by their components in $J_e=J_{11}$. Indeed, the components  $e'J'=J'_{21}\oplus J'_{22}$ are contained in both $\Ac$ and  $\Bc$. If the element
$x_{12}\in J_{12}$ belongs to  $\Ac$, then $\mu(x_{12}J')=0$. On the other hand,  as $W\in J_{11}$, we have $x_{12}W=0$. Hence $\mu(x_{12}J)=
\mu(x_{12}W)+\mu(x_{12}J')=0$, and $x_{12}\in \Bc$.

Since  $\Bc\subset\Ac$, the subgroup $ H(e)+\Bc$ is contained in the subgroup  $H(e)+\Ac$. Then
 \begin{multline}\label{xxxxx}  \Ind(\xi_{\theta,\mu}, H(e)+\Bc, G') =\\
  \Ind(\Ind(\xi_{\theta,\mu}, H(e)+\Bc, H(e)+\Ac), H(e)+\Ac, G').
                \end{multline}
Decompose the character
$\phi_{\theta,\mu'}=\Ind(\xi_{\theta,\mu}, H(e)+\Bc, H(e)+\Ac)$ into irreducible components.

Let  $\Fc$ stands for the set of all   $\bar{\nu}\in \Ac^*$ such that
$\bar{\nu}\vert_\Bc = \mu\vert_\Bc$. Observe that any  $\bar{\nu}\in\Fc$ annihilates each Pierce component of algebra  $\Ac$
except the first one (since this components belong to  $\Bc$, and $\mu\in J_e $ annihilates them). As $H(e)$ acts trivially on the first Pierce component, we verify  $\bar{\nu}(hy)=\bar{\nu}(yh)=\bar{\nu}(y)$ for each  $h\in H(e)$.

Let us show that  $\Ac^2\subseteq \Bc$. Indeed, if  $y_1,y_2\in \Ac$, then  $\mu'(y_1J')=\mu'(y_2J')=0$.
Hence $\mu(y_1y_2J)=\mu'(y_1(y_2J))=\mu'(y_1J')=0$. This proves   $\Ac^2\subseteq \Bc$.
Moreover, for any   $\bar{\nu}\in \Fc$, we obtain $\bar{\nu}(y_1y_2)=\mu(y_1y_2)=\mu'(y_1y_2)=0$.

The following formula
\begin{equation}\label{zzzzz}
 \xi_{\theta,\bar{\nu}}(h+y)= \theta(h)\varepsilon^{\bar{\nu}(y)},\quad h\in H(e),\quad y\in \Ac
\end{equation}
defines a linear character of the subgroup  $H(e)+\Ac$.
Really, let  $g_1=h_1+y_1$ and $g_2=h_2+y_2$, where $h_1,h_2\in H(e)$ and
$y_1,y_2\in A$. Then
\begin{multline*}
  \xi_{\theta,\bar{\nu}}(g_1g_2)= \xi_{\theta,\bar{\nu}}(h_1h_2+h_1y_2+y_1h_2+y_1y_2)=
\theta(h_1h_2)\eps^{\bar{\nu}(h_1y_2)}\eps^{\bar{\nu}(y_1h_2)}\eps^{\bar{\nu}(y_1y_2)}=\\
\theta(h_1)\theta(h_2)\eps^{\bar{\nu}(y_2)}\eps^{\bar{\nu}(y_1)}= \xi_{\theta,\bar{\nu}}(g_1) \xi_{\theta,\bar{\nu}}(g_2).
\end{multline*}
Any character  $ \xi_{\theta,\bar{\nu}}$ being restricted on  $H(e)+\Bc$ coincides with
$\xi_{\theta,\mu}$. By the  Frobenius reciprocity theorem, each   $ \xi_{\theta,\bar{\nu}}$,~ $\bar{\nu}\in \Fc$ is included in  $\phi_{\theta,\mu'}$.
The number of elements of $\Fc$ coincides with the degree of character  $\phi_{\theta,\mu'}$.
Therefore $$\phi_{\theta,\mu'} =\sum_{\bar{\nu}\in\Fc}  \xi_{\theta,\bar{\nu}}.$$
Substituting in (\ref{xxxxx}), we obtain
\begin{equation}\label{xiind}
    \Ind(\xi_{\theta,\mu}, H(e)+\Bc, G') = \sum_{\bar{\nu}\in\Fc}\Ind(\xi_{\theta,\bar{\nu}},  H(e)+\Ac, G').
   \end{equation}
Let us show that each summand of this sum is a sum of supercharacters.
 Since  $\Bc= J_{\mu,\rt}\cap J'\subset \Ac\subset J'$, for any  $\bar{\nu}\in\Fc$ there exists  $\nu\in J^*$ that is equal to
 $\bar{\nu}$ being restricted to  $\Ac$, and is equal to $\mu$ being restricted to $J_{\mu,\rt}$. As  $\mu $ and  $\bar{\nu}$ are equal on all Pierce components except the first one, we can choose   $\nu\in J_e$.
   The difference $\nu-\mu$   annihilates on $J_{\mu, \rt}$. The Lemma   \ref{lemann} implies that  $\nu-\mu\in J\mu$. Since $\nu,\mu\in J_e$, we have $$\nu-\mu=e\nu e - e\mu e= e(\nu-\mu)e\in eJ\mu e= eJe\mu.$$
 That is $\nu-\mu\in J_e\mu$, and  $\nu\in N_e\mu$, where $N_e=1+J_e$.
It is also true the opposite: if $\nu\in N_e\mu$, then  $\nu\vert_\Bc = \mu'\vert_\Bc$.

Since  $J'$ is a  $H\times H$- invariant ideal in  $J$, the group  $G$ acts on  $J'$ (and $J^{'*}$) by right ans left multiplications.
 The equality  $\nu=a\mu$, where $a\in N_e$, admits the restriction on  $J'^*$ in the form    $\nu'=a\mu'$.
 This implies that  $J'_{\nu',\rt} = J'_{\mu',\rt}=\Ac$. Therefore
 $$ \Ind(\xi_{\theta,\bar{\nu}}, H(e)+\Ac, G') = \Ind(\xi_{\theta,\nu'},  H(e)+J'_{\nu',\rt}, G').$$
According to Lemma  \ref{indind}, this character is a sum of supercharacters. This proves the statement in Item  2.1. \\
\emph{Item  2.2.} Suppose that  $W\subseteq J_{12}$.  Then $J_{11}=J'_{11}$,~  $J_{21}=J'_{21}$,~ $J_{22}=J'_{22}$, and
$J_{12} = W\oplus J'_{12}$.\\
As in the previous item, the restriction   $\Res\,\chi_\al$ is calculated by the formula  (\ref{restr}).
The subspace  $J_{\la,\rt}+J'$ is invariant with respect to  $H(e)\times H(e)$. Denote by $W_0$ the $H(e)\times H(e)$-invariant subspace in $ J_{12}$ such that
\begin{equation}\label{jwzero}
J = W_0\oplus (J_{\la,\rt}+J').
\end{equation}
  We consider that
$G'\backslash G/G_\al$ is represented by the elements  $s=1+z$, where $z\in W_0$.
Let us calculate   $sG_\al s^{-1}\cap G'$.

 Let  $g\in G_\al$. Then $g=h+x$, where $h\in H(e)$ and $x\in J_{\la,\rt}$.
 Let us show that   $sgs^{-1}$ belongs to  $G'$ if and only if  $h$ commutates with
$s$ (equivalently, with  $z$) and $sxs^{-1}\in  J'$.
 Indeed,  $s(h+x)s^{-1}= shs^{-1}+sx s^{-1}$.
Observe that since  $J_{12}^2=0$,  we have  $s^{-1}=1-z$.
 The element  $$shs^{-1}=(1+z)h(1-z)=h+[z,h]$$  belongs to  $h+W_0$.
The element $x$ is represented in the form  $x=x_{12}+y$, where  $x_{12}\in J_{12}$, and $y$ is a sum of all other Pierce components of  $x$.
By conditions of Item  2.2,  $y\in J'$. Since  $J'$ is an ideal in $J$, the element  $sys^{-1}$ belongs to  $J'$.
The subalgebra  $J_{\la,\rt}$  is graded with respect to the  Pierce decomposition. As $x\in J_{\la,\rt}$, we have  $x_{12}\in J_{\la,\rt}$.  Since  $J_{12}^2=0$, we calculate $sx_{12}s^{-1}=x_{12}$.
Then  $sxs^{-1}=sx_{12}s^{-1}+sys^{-1}=x_{12}+sys^{-1}\in J_{\la,\rt}+J'$.
Summing up, we conclude that   $sgs^{-1}$  can be written in the form  $sgs^{-1}=h+[z,h]+sxs^{-1}$, where $[z,h]$, $sxs^{-1}$ are
components in decomposition  (\ref{jwzero}). Therefore,  $sgs^{-1}\in G'$ if and only if
$[z,h]=0$ and $sxs^{-1}\in J'$.

Denote by  $F$ the set of all  $h\in H(e)$ that commutes with  $s$.
 Then $ sG_\al s^{-1}\cap G'$ consists of all elements of the form  $s(h+x)s^{-1}$, where $h\in F$, and  $sxs^{-1}$ belongs to  $ sJ_{\la,\rt}s^{-1}\cap J'$. That is  $$sG_\al s^{-1}\cap G'= F+(J_{s\la s^{-1},\rt}\cap J').$$

Let us show that  $ J_{s\la s^{-1},\rt}\cap J'=  J'_{s\la' s^{-1},\rt}$.
Obviously, the left hand side of this equality is contained in the right hand side. Let us show the opposite inclusion.
If $y\in J'_{s\la' s^{-1},\rt}$, then  $s\la s^{-1} (yJ')=0$.
Recall that  $J=J_{12}+J'$. To prove $y\in J_{s\la s^{-1},\rt}$ it is sufficient to show that $s\la s^{-1} (yJ_{12})=0$. Indeed, $s\la s^{-1} (yJ_{12})= J_{12} s\la s^{-1} (y)$. As $J_{12}^2=0$, we have $J_{12}s=J_{12}$. We obtain
$s\la s^{-1} (yJ_{12})=J_{12}\la s^{-1}(y)$. Finally, $J_{12}\la=eJe'\la=0$ since  $\la\in J^*_e$.

So,  $sG_\al s^{-1}\cap G' = F+ J'_{s\la' s^{-1},\rt}$.
The equality  (\ref{restr}) takes the form
\begin{equation}\label{casetwo}
  \Res\,\chi_\al = \sum_{\mu} \Ind(\xi_{\theta,\mu'}, F+J'_{\mu',\rt}, G'),
\end{equation}
where $\mu'$ runs through all  $s\la' s^{-1}$, ~ $s=1+z$,~ $z\in W_0$.

 It follows from Lemma \ref{indind} that each summand in the sum  (\ref{casetwo}) is a sum of supercharacters.  This proves the statement in Item  2.2.\\
 \emph{Item 2.3.} Let $W\subseteq J_{21}$ or  $W\subseteq J_{22}$.\\
 As $\la\in J_e^*$, we have $e'J\in J_{\la,\rt}$. Therefore  $W\subset J_{\la,\rt}$.
 In this case,  $J'+ J_{\la,\rt}=J$, and, hence, $G'\backslash G/G_\al=\{1\}$.

 The restriction  $\Res\,\chi_\al$ coincides with  $\Ind(\xi_{\theta,\la}, G_\al\cap G,G')$.
 As in Item  2.1,  $G_\al = H(e)+\Bc$, where $\Bc=J_{\la,\rt}\cap J'$.
 The subgroup  $G_\al$ is contained in  $H(e)+\Ac$, where $\Ac=J'_{\la',\rt}=\{y\in J': \la(yJ')=0\}$.

 Suppose that  $W\subseteq J_{22}$. Recall $\la\in J_e^*$, then $\la(yJ_{22})=0$ for any  $y\in J$.
 We conclude that $\Bc=\Ac$,  and  $\Res\,\chi_\al=\Ind(\xi_{\theta,\la'}, H(e)+J'_{\la',\rt}, G')$.
 According to Lemma \ref{lemann}, this character is a sum of supercharacters.

  Suppose that $W\subseteq J_{21}$. As in Item  2.2, since  $\Bc\subset\Ac$, the subgroup  $ H(e)+\Bc$ is contained in $H(e)+\Ac$, and the formula  (\ref{xxxxx}) holds.
Decompose the character
$\phi_{\theta,\la'}=\Ind(\xi_{\theta,\la'}, H(e)+\Bc, H(e)+\Ac)$ into irreducible components.

 Denote by  $\Fc$ the set of all  $\bar{\nu}\in \Ac^*$ such that  $\bar{\nu}|_\Bc= \bar{\la}|_\Bc$.
 Each element $\bar{\nu}\in \Fc$ we extend to   $\nu\in J^*$ that is equal to $\la$ being restricted to  $J_{\la,\rt}$. Then $\nu-\la$ is equal to zero on  $J_{\la,\rt}$. Applying Lemma  \ref{lemann}, we obtain  $ \nu-\la\in J\la$. Linear forms of   $J'\la$ annihilate both  $A$ and $J_{\la,\rt}$. Since  $J=W\oplus J'$, we  have  $\nu-\la\in W\la$.

 The set   $\Fc$ is a projection of elements $(1+W)\la$ on  $\Ac$. This argues that the group  $H(e)$ trivially act on  $\Fc$ by right multiplication and nontrivially by left multiplication  (this differs with Item  2.1). Consider the set of orbits with respect to the adjoint action of the group  $H(e)$  on  $\Fc$. Let $\{\bar{\nu}_i\}$ be the set of orbit representatives,  $\{H_i\}$ be their stabilizers in  $H(e)$, and $\theta_i$ be the restriction of $\theta $ on $H_i$. The formula  (\ref{zzzzz})  determines a linear character  $\xi_i=\xi_{\theta_i,\bar{\nu}_i}$ of the subgroup  $H_i+\Ac$. According the Mackey method
 of description of irreducible representations of semidirect products  \cite[Proposiiton 25]{Serr},  the induced characters
 $\chi_i=\Ind(\xi_i, H_i+\Ac, H(e)+\Ac)$ are irreducible and pairwise nonequivalent.

 Applying the  Frobenius reciprocity theorem, easy to verify that  $\chi_i$ includes in decomposition of
 $\phi_{\theta,\la'}$ into a sum of irreducible characters. Let us calculate the sum of degrees of these characters:
   $$\sum \mathrm{deg}\,\chi_i = \sum |H(e)/H_i|=|\Fc|= |\Ac/\Bc| = \mathrm{deg}\, \phi_{\theta,\la'}.$$
 From this we conclude  $$\phi_{\theta,\la'}=\sum \chi_i = \sum \Ind(\xi_i, H_i+\Ac, H(e)+\Ac).$$
Substituting in (\ref{xxxxx}), we obtain
$$\Res\,\chi_\al = \sum \Ind(\xi_i, H_i+\Ac, G).$$
It remains to show that the characters  $\Ind(\xi_i, H_i+\Ac, G)$ are sums of supercharacters.

Recall that    $\Ac=J'_{\la',\rt}$. Extend each  $\bar{\nu}_i\in \Fc$ to  $\nu_i\in J^*$ that equals to $\la$ on  $J_{\la,\rt}$. As we saw above,  $\nu_i=(1+w_i)\la$. Restricting to $J'$, we get $\nu'_i=(1+w_i)\la'$.
Hence, the right stabilizers of $\la'$ and $\nu'$ coincides and
$$\Ind(\xi_i, H_i+\Ac, G)=\Ind(\xi_{\theta_i,\nu_i'}, H_i+ J'_{\nu_i',\rt}, G').$$
Using Lemma  \ref{indind}, we conclude that this character is  a sum of supercharacters of the subgroup  $G'$. $\Box$\\
\Cor\Num. The product of supercharacters of  finite group of triangular type is a sum of supercharacters with nonnegative integer coefficients.\\
\Proof. Let  $G$ be a finite group of triangular type, and $\chi_1$, ~$\chi_2$ are its supercharacters. Then  $G\times G=H\times H + J\oplus J$ is also a finite group of triangular type, and $\chi_1\otimes\chi_2\in \mathrm{Fun}(G\times G)$ is its supercharacter.
The statement is verified using the restriction of   $\chi_1\otimes\chi_2$ on the diagonal subgroup $\{(g,g):~ g\in G\}$. $\Box$

\section{The Frobenius theorem for supercharacters}

  Let $G=H+J$ be a finite group of triangular type, $G'=H'+J'$ is its subgroup of triangular type. The scalar product on  $G$
  is defined as usual
 $$(f_1,f_2)=\frac{1}{|G|} \sum_{g\in G} f_1(g)\overline{f_2(g)}.$$
  Let $\phi$ be a superclass function  (i.e., a complex valued function constant on superclasses) on $G'$. Denote by  $\dot{\phi}$ the function on $G$ equal to  $\phi$ on $G'$ and zero out $G'$.

We define \emph{superinduction} as follows:
\ $$ \SInd\, \phi(g) = \frac{|H|}{|G|\cdot|G'|}\sum_{\tau\in \tG} \dot{\phi}(\rho(\tau)(g))=
\frac{|H|}{|G|\cdot|G'|}\sum_{a,b\in N,~ t\in H} \dot{\phi}(1+ta(g-1)bt^{-1}).$$
Easy to see that $\SInd\,\phi$ is a superclass function on $G$. \\
\Theorem\Num\label{last}. Let  $\psi$ be a superclass function on  $G$. Then   $$(\SInd\,\phi, \psi) = (\phi, \Res\,\psi).$$
\Proof. We obtain the proof by direct calculations:
\begin{multline*}
  (\SInd\,\phi, \psi) =\frac{|H|}{|G|^2\cdot|G'|}\sum_{\tau\in \tG,~ g\in G} \dot{\phi}(\rho(\tau)(g))\overline{\psi(g)} = \\
  \frac{|H|}{|G|^2\cdot|G'|}\sum_{\tau\in \tG,~ g\in G} \dot{\phi}(g)\overline{\psi(\rho(\tau)(g))}=
   \frac{|H|}{|G|^2\cdot|G'|}\cdot |N|^2\cdot|H|\sum_{\tau\in \tG,~ g\in G} \dot{\phi}(g)\overline{\psi(g)}=\\
   \frac{1}{|G'|} \sum_{g'\in G'} \phi(g')\overline{\psi(g')} = (\phi,\Res\,\psi). \Box
         \end{multline*}

   Let  $\{\chi_\al\}$ be the system of supercharacters of  finite group of triangular type $G=H+J$,
   and $\{\phi_\eta\}$ be the system of supercharacters of its subgroup of triangular type  $G'=H'+J'$. By Theorem  \ref{theoremrestr},  $$\Res\,\chi_\al=\sum_\eta m_{\al,\eta}\phi_\eta, ~ \quad  m_{\al,\eta}\in\Zb_+ .$$
    The system of supercharacters form a basis in the subspace of superclass functions  (see \cite{DI}).
   Since  $\SInd\,\phi_\eta$ is a superclass function on  $G$,  we obtain  $$ \SInd\,\phi_\eta=\sum_{\al}
 a_{\eta,\al}\chi_\al. $$
 The Theorem  \ref{last} implies  $$ a_{\eta,\al} = \frac{m_{\al,\eta}(\phi_\eta,\phi_\eta)}{( \chi_\al, \chi_\al)}.$$
 \Cor\Num. For any supercharacter  $\phi$ of subgroup of triangular type  $G'$, the superinduction  $\SInd\,\phi$ is a sum of supercharacters of the group  $G$ with nonnegative rational coefficients.

\end{document}